\documentclass[12pt]{article}

\usepackage{geometry}                
\usepackage{blindtext}               
\usepackage{graphicx}                
\usepackage[T1]{fontenc}             
\usepackage{tgpagella}               
\usepackage[symbol]{footmisc}        
\usepackage{hyperref}                
\usepackage{amsmath}                 
\usepackage{amssymb}
\usepackage{amsthm}                  
\usepackage{rotating}                
\usepackage{multirow}                
\usepackage{booktabs}                
\usepackage{adjustbox}               
\usepackage{mathtools}
\usepackage{xcolor}                  
\usepackage{colortbl}                
\usepackage{soul}
\usepackage[backend=biber,giveninits=false,uniquename=false, maxnames=99, uniquelist=false]{biblatex} 
\usepackage{subcaption}              
\usepackage{siunitx}                 
\usepackage{datetime}                
\usepackage{authblk}                 

\newcommand{\nf}[1]{\mathbf{#1}} 

\newcommand{\red}[1]{\textcolor{red}{#1}}

\newcommand{\Z}{\mathbb{Z}}
\newcommand{\Q}{\mathbb{Q}}
\newcommand{\R}{\mathbb{R}}
\newcommand{\C}{\mathbb{C}}
\newcommand{\F}{\mathbb{F}}

\DeclareMathOperator{\diag}{\mathrm{diag}}
\DeclareMathOperator{\Mat}{\mathrm{Mat}}

\DeclareMathOperator{\im}{\mathrm{Im}}

\newtheorem{theorem}{Theorem}
\newtheorem{corollary}[theorem]{Corollary}
\newtheorem{lemma}[theorem]{Lemma}
\newtheorem{proposition}[theorem]{Proposition}
\theoremstyle{definition}
\newtheorem{definition}[theorem]{Definition}
\theoremstyle{remark}
\newtheorem*{remark}{Remark}

\title{Determinantal Upper Bounds for Ternary Matrices}
\author[1,2]{Akihiro Munemasa}
\author[1,3]{Guillermo Núñez Ponasso}

\affil[1]{\footnotesize Graduate School of Information Sciences, Tohoku University, Sendai, Miyagi, Japan}
\affil[2]{\footnotesize
School of Science, China University of Geosciences, Beijing, China
}
\affil[3]{\footnotesize Dept. of Electrical \& Computer Engineering, Worcester Polytechnic Institute, Worcester, MA, USA}

\date{\today}

\begin{document}
\maketitle
\begin{abstract}
    We obtain a new, sharper upper bound for the determinant of a matrix of order $n\equiv 2\pmod 3$ with entries in the complex third-roots of unity. This bound is an analog of Ehlich's bound for binary matrices of order $n\equiv 3\pmod 4$ and it is obtained by first applying an appropriate normalization, or balancing, procedure for matrices with root-of-unity entries and then  generalizing Ehlich's argument. This new bound achieves an asymptotic fraction of approximately 74\% of the Hadamard bound.
\end{abstract}

\section{Introduction}

Denote by $\mu_{\ell}=\{1,\zeta_{\ell},\dots,\zeta_{\ell}^{\ell-1}\}$ the group of complex $\ell$-th roots of unity, where $\zeta_{\ell}$ denotes a primitive $\ell$-th root of unity. In the particular case where $\ell=3$ we denote $\omega=\zeta_3$. We fix once and for all the embedding $\omega=\exp(2\pi i/3)$ and note that all the results we present are independent of the choice of embedding.

The entrywise product, or Schur product, $u\circ v$ of two vectors $u$ and $v$ of the same length is the vector whose $i$-th entry is $(u\circ v)_i=u_iv_i$.

The identity matrix of order $n$ will be denoted $I_n$, and the all-ones matrix of order $n$ will be denoted $J_n$. The all-ones vector of length $n$ is denoted $\mathbf{1}_n$.

A monomial matrix is an element of $\Mat_{n}(\C)$ with exactly one non-zero element in each row and column. We denote by $\mathcal{M}_{n}(\ell)$ the set of monomial matrices whose non-zero entries are taken from the set $\mu_{\ell}$. We say that two matrices $M$ and $N$ in $\Mat_n(\mu_{\ell})$ are monomially-equivalent if and only if there exists a pair of matrices $P,Q\in\mathcal{M}_{n}(\ell)$ such that
\[N = PMQ^*.\]
Clearly, if $M$ and $N$ are monomially-equivalent then $|\det M|=|\det N|$.

Hadamard's determinant bound \cite{Hadamard93} states that the determinant of an order $n$ matrix $M$ with entries in the complex unit disk satisfies $|\det M|^2\leq n^n$, with equality in the bound if and only if $MM^*=nI_n$. In the same paper, Hadamard asked the question of determining the maximal determinant of a matrix with entries $\pm1$. The $\pm1$ matrices attaining the Hadamard bound with equality are known as Hadamard matrices, and they can only exist at orders $n=1,2$ or $n$ a multiple of $4$. Currently, the smallest open case for the existence of an Hadamard matrix is $n=2396$, \cite{EKLoC26}.
Sharper upper bounds have been found for $\pm 1$ matrices of order $n>2$ with $n\not\equiv 0\pmod 4$. Barba \parencite{Barba33} showed that a $\pm 1$ matrix $M$ of odd order $n$ satisfies $|\det M|^2\leq (2n-1)(n-1)^{n-1}$, with equality if and only if $NN^{\intercal}=(n-1)I_n+J_n$. \footnote{There exists a related determinantal upper bound for $n\equiv 2\pmod 4$; \cite{Ehlich64a, Wojtas64}}  Equality in Barba's bound implies in particular that $n\equiv 1\pmod 4$, hence this bound can be further refined for $n\equiv 3\pmod 4$. In \cite{Ehlich64b}, Ehlich analyzed the  maximal determinant of matrices in the finite set
\begin{equation}\label{eq:Cm}
	\mathcal{C}_m=\{C\in\Mat_{m}(3+4\Z): C=(c_{ij})_{ij},c_{ii}=n, C=C^* \text{ is positive definite}\},
\end{equation}
and maximized the determinant of a matrix in $\mathcal{C}_m$ to find a determinantal upper bound for matrices with entries $\pm 1$ of order $n\equiv 3\pmod{4}$.

The maximal determinant problem for $\pm1$ matrices has been generalized to matrices over the fourth roots of unity in \cite{Cohn96} and subsequently to matrices over arbitrary roots of unity \cite{Ponasso23,Ponasso25}. A generalization of Barba's bound, depending on the minimal absolute value of a sum of $n$ roots of unity, was shown in \cite{Ponasso23,Ponasso25}. For ternary matrices, i.e. matrices with entries in the set of third roots of unity, this generalized Barba bound is identical to the original, and it can only be met with equality if $n\equiv 1\pmod 3$, hence a sharpening like that of Ehlich could be possible for ternary matrices of order $n\equiv 2\pmod 3$. The main purpose of this paper is to obtain such an upper bound.

Here, we will proceed in a similar way to Ehlich's proof and with as much generality as we possible, first re-proving results in Ehlich's paper while at the same time deriving new results pertaining to matrices with entries in $\mu_3$, and finally analyzing the latter case in detail while referring to \cite{Ehlich64b} for the analogous calculations.

We will often use the following corollary of Sylvester's Criterion:

\begin{lemma}\label{lem:aux1} Let $A$ be a Hermitian matrix with $\det A>0$. If the submatrix of $A$ obtained by deleting the last row and column is positive definite, then $A$ is positive definite.
\end{lemma}
\begin{proof}
	Let $A'$ be the submatrix obtained from $A$ by deleting its last row and column. Since $A'$ is positive definite, all the principal minors of $A'$ are positive definite by \cite[Cor. 7.1.5]{HJ90}. This, together with the assumption that $\det A>0$, implies that every leading principal minor of $A$ is positive. By Sylvester's Criterion \cite[Thm. 7.2.5(b)]{HJ90} we conclude that $A$ is positive definite.
\end{proof}

Our main result is the following determinantal upper bound:
\begin{theorem}\label{thm:main}
	Let $M\in\Mat_n(\mu_3)$. If $n\geq 2$ is an integer with $n\equiv 2\pmod{3}$, then
	\[|\det M|^2\leq
	\begin{cases}
        D_n(2) & \text{if } n=2,\\
        D_n(3) & \text{if } 5\leq n \leq 23,\\
		D_n(4) & \text{if } n\geq 26,
	\end{cases}
	\]
	where 
	\[D_n(s)=(n-2)^{n-s}(n-2+3r)^u(n+1+3r)^v \left(1-\frac{ur}{n-2+3r}-\frac{v(r+1)}{n+1+3r}\right),\]
	$r=\lfloor n/s\rfloor$, $u=s(r+1)-n$, and $v=s-u$.
\end{theorem}

To prove this result, we will first introduce an appropriate normalization procedure for matrices with entries in the set of $\ell$-th roots of unity $\mu_{\ell}$. Normalization will ensure that Gram matrices belong to a set $\mathcal{K}_n$, for which we will be able to characterize its maximal determinant matrices as block matrices indexed by partitions of $n$. We will then determine the partition shapes associated to maximal determinant matrices and thus obtain a set of candidate maximal determinant matrices of size $n$. Finally, for the case $\ell=3$ we determine a unique maximal determinant of $\mathcal{K}_n$ for every $n$, leading to an upper bound for matrices with entries in $\mu_3$.
\section{Balanced matrices and their Gram matrices}
The aim of this section is to determine an appropriate set of matrices $\mathcal{G}$, which will contain all Gram matrices of matrices in $\Mat_n(\mu_3)$ up to equivalence.

Namely, we show that up to monomial equivalence, any matrix $M\in \Mat_n(\mu_3)$ has a Gram matrix $G=MM^*\in\mathcal{G}$ that satisfies certain congruence relations. We describe this monomial transformation by introducing the concept of balanced vectors in $\mu_{\ell}^n$ and balanced matrices $\Mat_n(\mu_{\ell} )$.

\begin{definition}
	Let $\ell>1$ and $n\geq 1$ be integers. A vector $v\in(\mu_{\ell})^n$ written as $v=[\zeta_{\ell}^{e_1},\zeta_{\ell}^{e_2},\dots,\zeta_{\ell}^{e_n}]$ is said to be balanced if $\sum_{i=1}^n e_i\equiv 0\pmod{\ell}$.
\end{definition}

For $v\in\mu_{\ell}^n$ and $i\in\{0,1,\dots,\ell\}$ let $m_i(v)$ be the multiplicity of $\zeta_{\ell}^i$ in $v$, i.e. the number of times the entry $\zeta_{\ell}^i$ appears in $v$. Denote $m(v):=(m_0(v),m_1(v),\dots,m_{\ell-1}(v))$.

A vector $v\in\mu_3^n$ with $(a,b,c):=m(v)$ is balanced if and only if $0\cdot a+1\cdot b+2\cdot c\equiv 0\pmod{3}$, which is equivalent to $b\equiv c\pmod{3}$. Note that if $n\equiv 0\pmod{3}$, then the condition $b\equiv c\pmod{3}$ combined with $a+b+c=n$ implies $a\equiv b\equiv c\pmod{3}$. Over $\mu_3$, the notion of balanced vectors introduced here coincides with the one in \cite{Ponasso25}.

When $n$ is odd, a vector $v\in\mu_2^n$ with  $(a,b):=m(v)$ is balanced if and only if $0\cdot a+1\cdot b \equiv 0\pmod{2}$. So $v$ is balanced if and only if the number of $-1$s in $v$ is even. Since $n$ is odd, multiplying by $v$ by $-1$ the parities of $+1$s and $-1$s are exchanged. Therefore, the notion of balanced vectors for $\ell=2$ is equivalent to the notion of ``parity-normalized vector'' in \cite{Orrick05}, which instead requires the number of $+1$s to be even.

It is straightforward to prove the following lemma.
\begin{lemma}\label{lem:l1}{\ }
	\begin{itemize}
		\item[(i)] The all-ones vector $\mathbf{1}_n$ is a balanced vector in $\mu_{\ell}^n$ for all $\ell$.
		\item[(ii)] The complex conjugate $\overline{v}$ of a balanced vector $v\in\mu_{\ell}^n$ is balanced. More generally, the entrywise image $v^{\sigma}$ of $v$ under a field automorphism $\sigma$ of $\Q[\zeta_{\ell}]$ is balanced.
		\item[(iii)] The entrywise product $u\circ v$ of two balanced vectors $u$ and $v$ in $\mu_{\ell}^n$ is balanced.
	\end{itemize}
\end{lemma}

\begin{lemma} \label{lem:l2} Let $\ell>1$ and $n\geq 1$ be integers. Suppose that $\gcd(n,\ell)=1$ and $v\in\mu_{\ell}^n$, then there is a unique phase $\zeta_{\ell}^k\in\mu_{\ell}$ such that $\zeta_{\ell}^kv$ is a balanced vector.
\end{lemma}
\begin{proof} If $v = [\zeta_{\ell}^{e_i}]_i$, then the $i$-th entry of  $\zeta_{\ell}^kv$ is $\zeta_{\ell}^{e_i+k}$. Hence, $\zeta_{\ell}^{k}v$ is balanced if and only if
	\begin{equation}\label{eq:l2-1}
		\sum_{i=1}^n(k+e_i)=nk+\sum_{i=1}^{n}e_i\equiv 0\pmod{\ell}.
	\end{equation}
	Since $\gcd(n,\ell)=1$, there is an $a\in\{1,2,\dots,\ell-1\}$ such that $an\equiv 1\pmod{\ell}$. Therefore, rearranging (\ref{eq:l2-1}) and multiplying both sides by $a$ we find that choosing $k\in\Z$ such that
	\[k\equiv -a\sum_{i=1}^n e_i\pmod{\ell},\]
	then, by construction, $\zeta_{\ell}^kv$ is balanced. If both $\zeta_{\ell}^k v$ and $\zeta_{\ell}^{k'}v$ are balanced, then $nk+\sum_{i}e_i\equiv 0\equiv nk'+\sum_{i}e_i\pmod{\ell}$, hence $nk\equiv nk'\pmod{\ell}$ and multiplying by $a$ we find $k\equiv k'\pmod{\ell}$, so $\zeta_{\ell}^k=\zeta_{\ell}^{k'}$.
\end{proof}

By the above, if $\gcd(n,\ell)=1$, then for any  $M\in\Mat_n(\mu_{\ell})$ there is a unique diagonal matrix $P\in\mathcal{M}_{n}(\ell)$ such that $PM$ has all rows balanced. Furthermore, $|\det(PM)|=|\det M|$. These observations are sufficient to proceed with our proof of Theorem \ref{thm:main}. However, we can impose a stronger balancing condition for matrices as follows

\begin{definition} A matrix $M\in\Mat_n(\mu_{\ell})$ is called balanced if all of its rows and columns are balanced.
\end{definition}

\begin{lemma}\label{lem:l3}
	Let $\ell>1$ and $n\geq 1$ be integers. If $\gcd(n,\ell)=1$ and $M\in\Mat_n(\mu_{\ell})$ then there is a pair $(P,Q)$ of diagonal matrices in $\mathcal{M}_{n}(\ell)$ such that $PMQ\in\Mat_n(\mu_{\ell})$ is balanced.
\end{lemma}
\begin{proof}
	By Lemma~\ref{lem:l2} there is a unique diagonal matrix $P\in\mathcal{M}_{n}(\ell)$  such that $PM$ has all rows balanced. Similarly, there exists a unique diagonal matrix $Q\in\mathcal{M}_{n}(\ell)$ such that $PMQ$ has all columns balanced. It suffices to show that $PMQ$ has balanced rows.

	Write $Q=\diag(v)$ where $v\in\mu_{\ell}^n$, and let $u_r$ be the $r$-th row of $PM$. Then, the $r$-th row of $PMQ$ is $v\circ u_r$. By part (iii) of Lemma~\ref{lem:l1} it suffices to show that $v$ is balanced. Write $(PM)_{ij}=\zeta_{\ell}^{a_{ij}}$, and write $v=[\zeta_{\ell}^{e_i}]_i$. Since the $j$-th column of $PMQ$ is balanced, we have $\sum_{i=1}^n (a_{ij}+e_j)\equiv 0\pmod{\ell}$. Then, as in the proof of Lemma~\ref{lem:l2}, there is a unique choice for $e_j$ modulo $\ell$, namely
	\[e_j\equiv -a\sum_{i=1}^na_{ij}\pmod{\ell},\]
	where $an\equiv 1\pmod{\ell}$. Since $PM$ has all rows balanced, we have $\sum_{i,j=1}^n a_{ij}\equiv 0\pmod{\ell}$, hence
	\[\sum_{j=1}^ne_j\equiv -a\sum_{j=1}^n\sum_{i=1}^na_{ij}\equiv 0\pmod{\ell}.\]
	Therefore, $v$ is balanced and this concludes the proof.
\end{proof}
One can show furthermore that any other pair $(P',Q')$ of diagonal matrices in $\Mat_n(\mu_3)$ making $P'MQ'$ balanced satisfies $(P',Q')=(P\zeta^i,Q\zeta^{-i})$ for some phase $\zeta^i\in\mu_{\ell}$. This is shown most easily by writing matrices in ``logarithmic form'' as elements of $\Mat_n(\Z/\ell\Z)$.

We return to the case $\ell=3$. The following Lemma~determines arithmetic conditions that the entries of the Gram matrix of a balanced matrix in $\Mat_n(\mu_3)$ must satisfy.

\begin{lemma} \label{lem:l4} Let $u$ and $v$ be balanced vectors in $\mu_3^n$. Then, $u\cdot v\in n+3\Z[\omega]$.
\end{lemma}
\begin{proof}
	The inner product or $u$ and $v$ can be written as $u\cdot v=\mathbf{1}_n\cdot (\overline{u}\circ v)$, and by parts (ii) and (iii) of Lemma~\ref{lem:l1} the vector $\overline{u}\circ v$ is balanced. Therefore, without loss of generality, we may assume that $u=\mathbf{1}_n$. In this case, the inner product of $u$ and $v$ is equal to $\sum_{i=1}^n v_i=\sum_{i=1}^n \omega^{e_i}$. Denote $(a_0,a_1,a_2):=m(v)$. Since $v$ is balanced, we have $a_1\equiv a_2\pmod{3}$, which implies
	\begin{equation*}u\cdot v=a_0+a_1\omega+a_2\omega^2=n-(a_1-a_2)-3a_2+(a_1-a_2)\omega\in n+3\Z[\omega].\qedhere
	\end{equation*}
\end{proof}

\begin{lemma}\label{lem:l4-1} If $u\in 2+3\Z[\omega]$ and $|u|\leq 2$, then $u\in\{-1,2\}$.
\end{lemma}
\begin{proof}
	Write $u=2+3(x+y\omega)$, then
	\[4\geq |u|^2\geq\im(u)=\im(3y\omega)^2=\frac{27}{4}y^2.\]
	Therefore $|y|<1$, which implies $y=0$. So, $u\in 2+3\Z$ and since $|u|\leq 2$ it follows $u\in\{-1,2\}$.
\end{proof}

\begin{corollary}\label{cor:c1} Let $n\equiv 2\pmod{3}$. If $u,v\in\mu_3^n$ are balanced and 
	$|u\cdot v|\leq 2$, then 
	$u\cdot v\in\{-1,2\}$.
\end{corollary}
\begin{proof} This is a consequence of Lemmas~\ref{lem:l4} and~\ref{lem:l4-1}.
\end{proof}

Similarly, we obtain the following fact used in Ehlich's arguments.
\begin{lemma}\label{lem:l5} 
	Let $u,v\in\mu_2^n$ be balanced vectors, then $u\cdot v\in n+4\Z$. In particular, if $n\equiv 3\pmod{4}$ and $|u\cdot v|\leq 3$, then $u\cdot v\in\{-1,3\}$.
\end{lemma}
\begin{proof}
	As in the proof of Lemma~\ref{lem:l4}, we may assume without loss of generality that $u=\mathbf{1}_n$. If $(a,b):=m(v)$ then  $u\cdot v=a-b$, where $b\equiv 0\pmod{2}$ since $v$ is balanced. Then
	$u\cdot v
	\equiv a+b=n \pmod{4}$.
	If $n\equiv 3\pmod{4}$, we have $u\cdot v\in 3+4\Z$. Finally, the only elements in $3+4\Z$ of absolute value $\leq 3$ are precisely $-1$ and $3$.
\end{proof}

Denote
\begin{equation}
	\mathcal{P}_m := \{C\in\Mat_m(\C):C=(c_{ij})_{ij}, C=C^* \text{ is positive definite, }c_{ii}=n\}.
\end{equation}
We conclude this section with the determination of a suitable set containing all of the Gram matrices relevant to our study.
\begin{proposition}\label{prop:1}
	Let $m$ and $n$ be integers with $n\equiv 2\pmod{3}$ and $m\leq n$. If $M\in\Mat_{m\times n}(\mu_3)$ is full rank, then there exists a 
	diagonal matrix $P$ in $\mathcal{M}_n(3)$ such that $(PM)(PM)^*\in \mathcal{G}_m$, where
	\begin{equation}\label{eq:Gm}
		\mathcal{G}_m:=\Mat_m(2+3\Z[\omega])\cap\mathcal{P}_m.
	\end{equation}
	In particular,
	\begin{equation}\label{eq:t1-1}
		\max_{M\in \Mat_{n}(\mu_3)}|\det M|^2\leq \max_{C\in\mathcal{G}_n} |\det C|.
	\end{equation}
\end{proposition}
\begin{proof}
	By Lemma~\ref{lem:l2} there is a matrix $P$ such that $PM$ has balanced rows. By Lemma~\ref{lem:l4}, the inner product of two rows of $PM$ is in the set $2+3\Z[\omega]$. Therefore, $(PM)(PM)^*\in\Mat_m(2+3\Z[\omega])\cap\mathcal{P}_m=\mathcal{G}_m$.
	If $M\in \Mat_{n}$, then 
	$(PM)(PM)^*\in\mathcal{G}_n$.
	Since $|\det M|^2=|\det(MM^*)|=|\det(PM)(PM)^*|$, the inequality (\ref{eq:t1-1}) follows.
\end{proof}

\begin{remark} Note that since the matrices in $\mathcal{G}_m$ are Hermitian and positive definite,  the determinant of any principal $2\times 2$ submatrix is
	\[\det\begin{bmatrix}
		n & c_{ij}\\
		\overline{c_{ij}} & n
	\end{bmatrix}=n^2 -|c_{ij}|^2>0,\]
	which implies that $n>|c_{ij}|$. Therefore, the set $\mathcal{G}_m$ is finite for all $m\leq n$.
\end{remark}

\section{Proof of the determinantal upper bound}

We begin with a generalization of Satz 2.1 in \cite{Ehlich64b}. This is an auxiliary result to show that maximal determinant matrices in the set $\mathcal{C}_m$ (\ref{eq:Cm}) must have entries belonging to the set $\{-1,3\}$. With our generalization, we will be able to recover this result and show analogously that maximal determinant matrices in $\mathcal{G}_m$ have entries belonging to $\{-1,2\}$.

Fix integers $\ell$ and $n$ with $n>\ell>1$, and let $R$ be a discrete subset of $\C$ satisfying
\begin{equation}\label{eq:norm}
	\{z\in R: |z|\leq \ell\}=\{-1,\ell\}.
\end{equation}
We denote
\[\mathcal{K}_m = 
\Mat_m(R)\cap\mathcal{P}_m.
\]
Notice that for $m\leq n$, $\mathcal{K}_m$ is a non-empty finite set, because $(n+1)I_m-J_m\in\mathcal{K}_m$, and positive definiteness implies that all off diagonal entries of matrices in $\mathcal{K}_m$ are of modulus less than $n$.
For example, for $(\ell, R)=(2,2+3\Z[\omega])$, and $(\ell,R)=(3,3+4\Z)$ we have by Lemmas~\ref{lem:l4-1} and~\ref{lem:l5} that \eqref{eq:norm} is satisfied.
\begin{lemma}\label{lem:S21} Let $m\leq n$ be a positive integer, $\Delta_m = \max\{\det C : C\in\mathcal{K}_m\}$ and $\widehat{C}_m\in\mathcal{K}_m$, with $\det \widehat{C}_m=\Delta_m$. If $D_m$ is the matrix obtained from $\widehat{C}_m$ by replacing its $(m,m)$-th entry by $\ell$, then
	\begin{align}
		\det D_m &\geq (n-\ell)^{m-2}(\ell n-1) \quad\text{ for } 2\leq m\leq n, \label{eq:Dm}\\
		\Delta_{m+1}-(n-\ell)\Delta_{m}&\geq\begin{cases}
			\ell n-1 & \text{ if } m=1,\\
			(n-\ell)\det D_m & \text{ if } 2\leq m\leq n-1.
		\end{cases}\label{eq:Delm}
	\end{align}
	Furthermore, $D_m$ is positive definite for $1\leq m\leq n$.
\end{lemma}
\begin{proof} We proceed by induction on $m$, $1\leq m\leq n$. Let $m=1$, then $\mathcal{K}_m=\{[n]\}$, hence $\widehat{C}_1=[n]$, $\Delta_1=n$, and $D_1=[\ell]$ is positive definite. By definition
	\[\mathcal{K}_2=\left\{\begin{bmatrix}
		n & c\\
		\overline{c} & n
	\end{bmatrix}: c\in R, |c|<n\right\}.\]
	So, we have
	\begin{equation*}
		\begin{split}
			\Delta_2 &= \max\{n^2-|c|^2: c\in R,|c|<n\}\\
			&=n^2-\min\{|c|^2:c\in R,|c|<n\}\\
			&=n^2-1.
		\end{split}
	\end{equation*}
	Hence,
	\begin{equation*}
		\begin{split}
			\widehat{C}_2&=\begin{bmatrix}n&-1\\-1&n\end{bmatrix}\quad\text{and,}\\
				D_2 &=\begin{bmatrix}n & -1\\ -1 &\ell\end{bmatrix}.
		\end{split}
	\end{equation*}
	Since $\det D_2=\ell n-1>0$, (\ref{eq:Dm}) holds and $D_2$ is positive definite. Also,
	\begin{equation*}
		\begin{split}
			\Delta_2-(n-\ell)\Delta_1&=\det \widehat{C}_2-(n-\ell)\det\widehat{C}_1\\
			&=n^2-1 - (n-\ell)n \\&= \ell n-1,
		\end{split}
	\end{equation*}
	which shows (\ref{eq:Delm}) for $m=1$. 

	Now suppose $m\geq 2$. 
	By the inductive hypothesis
	\begin{equation*}
		\begin{split}
			\Delta_m-(n-\ell)\Delta_{m-1}&\geq
			\begin{cases}
				\ell n-1 &\text{if } m=2,\\
				(n-\ell)\det D_{m-1} & \text{if } 3\leq m\leq n
			\end{cases}\\
			&\geq \begin{cases}
				\ell n-1 &\text{if } m=2,\\
				(n-\ell)^{m-2}(\ell n-1) & \text{if } 3\leq m\leq n
			\end{cases}\\
			&=(n-\ell)^{m-2}(\ell n -1).
		\end{split}
	\end{equation*}
	Write
	\[\widehat{C}_m=\begin{bmatrix}
		C' & \mathbf{c}\\
		\mathbf{c}^* & n
	\end{bmatrix},
	\]
	then
	\[D_m=\begin{bmatrix}
		C' & \mathbf{c}\\
		\mathbf{c}^* & \ell
	\end{bmatrix}.\]
	Since $\widehat{C}_m\in\mathcal{K}_m$ is positive definite, its leading principal submatrix $C'$ must also be positive definite, hence $C'\in\mathcal{K}_{m-1}$. This implies that $\det C'\leq \Delta_{m-1}$. Then,
	\begin{equation*}
		\begin{split}
			\det D_m &= 
			\det \begin{bmatrix}
				C' & \mathbf{c}\\
				\mathbf{c}^* & \ell
			\end{bmatrix}
			+\det
			\begin{bmatrix}
				C' & \mathbf{c}\\
				\mathbf{0} & n-\ell
			\end{bmatrix}
			-(n-\ell)\det C'\\
			&=\det \begin{bmatrix}
				C' & \mathbf{c}\\
				\mathbf{c}^* & n
			\end{bmatrix}-(n-\ell)\det C'\\
			&=\det\widehat{C}_m - (n-\ell)\det C'\\
			&=\Delta_m - (n-\ell)\Delta_{m-1}\\
			&\geq (n-\ell)^{m-2}(\ell n-1).
		\end{split}
	\end{equation*}
	This proves (\ref{eq:Dm}), and in particular $\det D_m>0$, therefore, by Lemma~\ref{lem:aux1} $D_m$ is positive definite.

	Define $C\in \Mat_{m+1}(R)$ by
	\[
		C = \begin{bmatrix}
			C' & \mathbf{c} & \mathbf{c}\\
			\mathbf{c}^* & n & \ell\\
			\mathbf{c}^* & \ell & n
		\end{bmatrix}.
	\]
	Then,
	\begin{equation*}
		\begin{split}
			\det C &=\det
			\begin{bmatrix}
				C' & \mathbf{c} & \mathbf{c}\\
				\mathbf{c}^* & n & \ell\\
				\mathbf{0} & 0 & n-\ell
			\end{bmatrix}
			+
			\det
			\begin{bmatrix}
				C' & \mathbf{c} & \mathbf{c}\\
				\mathbf{c}^* & n & \ell\\
				\mathbf{c}^* & \ell & \ell
			\end{bmatrix}\\
			&=
			(n-\ell)\det
			\begin{bmatrix}
				C' & \mathbf{c} \\
				\mathbf{c}^* & n \\
			\end{bmatrix}
			+
			\det
			\begin{bmatrix}
				C' & \mathbf{c} & \mathbf{c}\\
				\mathbf{c}^* & n & \ell\\
				\mathbf{0} & \ell - n & 0
			\end{bmatrix}\\
			&= (n-\ell)\det\widehat{C}_m + (n-\ell) \det
			\begin{bmatrix}
				C' & \mathbf{c}\\
				\mathbf{c}^* & \ell\\
			\end{bmatrix}\\
			&= (n-\ell)\Delta_m + (n-\ell)\det D_m\\
			&>0.
		\end{split}
	\end{equation*}
	Since the submatrix $\widehat{C}_m$ of $C$ is positive definite, Lemma \ref{lem:aux1} implies that $C$ is also positive definite, hence $C\in\mathcal{K}_{m+1}$. Then, by the above,
	\begin{equation*}
		\begin{split}
			\Delta_{m+1} - (n-\ell)\Delta_{m} &\geq \det C - (n-\ell)\Delta_m\\
			&= (n-\ell)\det D_m.
		\end{split}
	\end{equation*}
	This proves (\ref{eq:Delm}).
\end{proof}

\begin{corollary}[{\cite[Satz 2.1.]{Ehlich64b}}] Suppose $n\equiv 3\pmod{4}$, $1\leq m<n$, and let $g_m:=\max\{\det C:C\in\mathcal{C}_m\}$, with $\mathcal{C}_m$ given in \eqref{eq:Cm}. Then, $g_{m+1}>(n-3)g_{m}$.
\end{corollary}
\begin{proof}
	For $R=3+4\Z$ and $\ell=3$ we have $n\in R$, and hence $\mathcal{K}_m=\mathcal{C}_m$. Applying Lemma~\ref{lem:S21}, the result follows. 
\end{proof}
\begin{corollary} Suppose that $n\equiv 2\pmod{3}$, $1\leq m< n$, and let $g'_m:=\max\{\det C:C\in\mathcal{G}_m\}$, with $\mathcal{G}_m$  given in \eqref{eq:Gm}. Then $g'_{m+1}>(n-2)g'_{m}$.
\end{corollary}
\begin{proof}
	For $R=2+3\Z[\omega]$ and $\ell=2$ we have $n\in R$, and hence $\mathcal{K}_m=\mathcal{G}_m$. Applying Lemma~\ref{lem:S21}, the result follows.
\end{proof}

In the definition of $D_m$ above, instead of replacing the $(m,m)$-th entry with $\ell$ we could have replaced it with any other positive real element $\alpha$ of $R$ and obtained 
the inequality 
$\Delta_m-(n-\alpha)\Delta_{m-1}\geq (n-\alpha)\det D_m$
which is similar to \eqref{eq:Delm}.
However, since $\ell$ is the smallest positive real element of $R$ this inequality is weaker than the one in Lemma~\ref{lem:S21}. If we replaced by a negative element of $R$ such as $-1$ instead of $\ell$, we would find that at the base case the corresponding matrix
\[D_2=\begin{bmatrix}
	n & -1\\
	-1 & -1
\end{bmatrix}
\]
has determinant $-n-1<0$, and we would not obtain an analogous result for all $m$. Finally, if we had chosen a complex element of $R$ in the definition of $D_m$, then $D_m$ would not be Hermitian.

\begin{lemma} \label{lem:aux2} Consider the partitioned Hermitian matrix
	\[M=\begin{bmatrix}
		A & B\\
		B^* & C
	\end{bmatrix}
	\in \Mat_{p+q}(\C)
	\quad \text{where }
	A\in\Mat_p(\C), 
	B\in\Mat_q(\C).
	\]
	If $\det C=0$ and the submatrix of $M$ obtained by deleting its last row and column is positive definite, then $\det M\leq 0$.
\end{lemma}
\begin{proof}
	If $\det M>0$, then every leading principal minor of $M$ is positive. By Sylvester's Criterion \cite[Cor. 7.1.5. and Thm. 7.2.5]{HJ90}, it follows that $M$ is positive definite. But then, this forces $\det C>0$ by \cite[Cor. 7.1.5]{HJ90}, contrary to our assumption.
\end{proof}

\begin{lemma}\label{lem:S22} Let $m$ be a positive integer with $m\leq n$. Suppose that $\widehat{C}=(c_{ij})\in\mathcal{K}_m$ satisfies $\det \widehat{C} =\Delta_m$. Then $c_{ij}\in\{-1,\ell\}$ for all $i,j\in\{1,\dots,m\}$ with $i\neq j$.
\end{lemma}
\begin{proof} By way of contradiction, suppose that $c_{ij}\not\in\{-1,\ell\}$, for some $i,j\in\{1,\dots,m\}$ with $i\neq j$. Then, up to conjugation by a permutation matrix we may assume that
	\[\widehat{C}=
	\begin{bmatrix}
		A & \mathbf{a}_1 & \mathbf{a}_2\\
		\mathbf{a}^*_1 & n & \alpha\\
		\mathbf{a}^*_2 & \overline{\alpha} & n
	\end{bmatrix},
	\]
	for some $\alpha\in R\setminus\{-1,\ell\}$. Then, by \eqref{eq:norm},
	\begin{equation}\label{eq:8a1}
		|\alpha|>\ell.
	\end{equation}
	We may further assume, swapping the last two rows and columns if necessary, that
	\begin{equation}\label{eq:8b1}
		\det\begin{bmatrix}
			A & \mathbf{a}_2\\
			\mathbf{a}_2^* & n
		\end{bmatrix}
		\geq
		\det\begin{bmatrix}
			A & \mathbf{a}_1\\
			\mathbf{a}_1^* & n
		\end{bmatrix}.
	\end{equation}
	Exchanging the last two rows and columns of $\widehat{C}$, the last statement of Lemma~\ref{lem:S21} implies that
	\begin{equation}\label{eq:8a3n}
		\begin{bmatrix}
			A & \mathbf{a}_2 & \mathbf{a}_1\\
			\mathbf{a}_2^* & n & \overline{\alpha}\\
			\mathbf{a}_1^* &\alpha & \ell
		\end{bmatrix}
	\end{equation}
	is positive definite. 
	Suppose
	\begin{equation*}
		\det 
		\begin{bmatrix}
			A & \mathbf{a}_1 & \mathbf{a}_2\\
			\mathbf{a}^*_1 & \ell & \alpha\\
			\mathbf{a}^*_2 & \overline{\alpha} & |\alpha|^2/\ell
		\end{bmatrix}
		> 0.
	\end{equation*}
	Since
	\begin{equation*}
		\begin{bmatrix}
			A & \mathbf{a}_1\\
			\mathbf{a}_1^* & \ell
		\end{bmatrix}
	\end{equation*}
	is a principal submatrix of the positive definite matrix \eqref{eq:8a3n}, it is
	positive definite. 
	By Lemma~\ref{lem:aux1} 
	the entire matrix above is positive definite. 
	But since the $2\times 2$ minor in the bottom right vanishes, this is a contradition. Hence,
	\begin{equation}\label{eq:8a3}
		\det 
		\begin{bmatrix}
			A & \mathbf{a}_1 & \mathbf{a}_2\\
			\mathbf{a}^*_1 & \ell & \alpha\\
			\mathbf{a}^*_2 & \overline{\alpha} & |\alpha|^2/\ell
		\end{bmatrix}
		\leq 0.
	\end{equation}
	Furthermore
	\begin{align}
		\det\begin{bmatrix}
			A & \mathbf{a}_2\\
			\mathbf{a}_2^* & \ell
		\end{bmatrix}
		&=\det
		\begin{bmatrix}
			A & \mathbf{a}_2\\
			\mathbf{a}_2^* & n
		\end{bmatrix}
		-\det
		\begin{bmatrix}
			A & \mathbf{a}_2\\
			\mathbf{0} & n-\ell
		\end{bmatrix}
		\notag\\
		&\geq
		\det
		\begin{bmatrix}
			A & \mathbf{a}_1\\
			\mathbf{a}_1^* & n
		\end{bmatrix}
		-\det
		\begin{bmatrix}
			A & \mathbf{a}_1\\
			\mathbf{0} & n-\ell
		\end{bmatrix}\quad \text{by (\ref{eq:8b1})}
		\label{eq:8a4}\\
		&=\det\begin{bmatrix}
			A & \mathbf{a}_1\\
			\mathbf{a}_1^* & \ell
		\end{bmatrix}
		>0.
		\label{eq:8a2}
	\end{align}
	Let
	\begin{equation*}
		C=
		\begin{bmatrix}
			A & \mathbf{a}_2 & \mathbf{a}_2\\
			\mathbf{a}^*_2 & n & \ell\\
			\mathbf{a}^*_2 & \ell & n
		\end{bmatrix}.
	\end{equation*}
	Then,
	\begin{align*}
		\det C
		&=
		\det
		\begin{bmatrix}
			A & \mathbf{a}_2 & \mathbf{a}_2\\
			\mathbf{a}^*_2 & n & \ell\\
			\mathbf{0} & 0 & n-\ell
		\end{bmatrix}
		+
		\det \begin{bmatrix}
			A & \mathbf{a}_2 & \mathbf{a}_2\\
			\mathbf{a}^*_2 & n & \ell\\
			\mathbf{a}^*_2 & \ell & \ell
		\end{bmatrix}\\
		&=
		(n-\ell)
		\det
		\begin{bmatrix}
			A & \mathbf{a}_2 \\
			\mathbf{a}^*_2 & n 
		\end{bmatrix}
		+\det
		\begin{bmatrix}
			A & \mathbf{a}_2 & \mathbf{a}_2\\
			\mathbf{0} & n-\ell & 0\\
			\mathbf{a}^*_2 & \ell & \ell
		\end{bmatrix}
		\displaybreak[0]\\
		&=(n-\ell)\det
		\begin{bmatrix}
			A & \mathbf{a}_2 \\
			\mathbf{a}^*_2 & n 
		\end{bmatrix}
		+(n-\ell)
		\det
		\begin{bmatrix}
			A & \mathbf{a}_2 \\
			\mathbf{a}^*_2 & \ell 
		\end{bmatrix}\\
		&\geq
		(n-\ell)\det
		\begin{bmatrix}
			A & \mathbf{a}_2 \\
			\mathbf{a}^*_2 & n 
		\end{bmatrix}
		+(n-\ell)
		\det
		\begin{bmatrix}
			A & \mathbf{a}_1 \\
			\mathbf{a}^*_1 & \ell 
		\end{bmatrix}
		\quad \text{ by (\ref{eq:8a4})}\\
		&>
		(n-\ell)\det
		\begin{bmatrix}
			A & \mathbf{a}_2 \\
			\mathbf{a}^*_2 & n 
		\end{bmatrix}
		+(n-|\alpha|^2/\ell)
		\det
		\begin{bmatrix}
			A & \mathbf{a}_1 \\
			\mathbf{a}^*_1 & \ell 
		\end{bmatrix}
		\quad\text{ by (\ref{eq:8a1}) and (\ref{eq:8a2})}\\
		&\geq
		(n-\ell)\det
		\begin{bmatrix}
			A & \mathbf{a}_2 \\
			\mathbf{a}^*_2 & n 
		\end{bmatrix}
		+(n-|\alpha|^2/\ell)
		\det
		\begin{bmatrix}
			A & \mathbf{a}_1 \\
			\mathbf{a}^*_1 & \ell 
		\end{bmatrix}
		+\det
		\begin{bmatrix}
			A & \mathbf{a}_1 & \mathbf{a}_2\\
			\mathbf{a}_1^* & \ell & \alpha\\
			\mathbf{a}_2^* & \overline{\alpha} & |\alpha|^2/\ell
		\end{bmatrix}
		\quad\text{ by (\ref{eq:8a3})}
		\displaybreak[0]\\
		&=
		(n-\ell)\det
		\begin{bmatrix}
			A & \mathbf{a}_2 \\
			\mathbf{a}^*_2 & n 
		\end{bmatrix}
		+\det
		\begin{bmatrix}
			A & \mathbf{a}_1 & \mathbf{a}_2\\
			\mathbf{a}_1^* & \ell & \alpha\\
			\mathbf{0} & 0 & n-|\alpha|^2/\ell
		\end{bmatrix}
		+\det
		\begin{bmatrix}
			A & \mathbf{a}_1 & \mathbf{a}_2\\
			\mathbf{a}_1^* & \ell & \alpha\\
			\mathbf{a}_2^* & \overline{\alpha} & |\alpha|^2/\ell
		\end{bmatrix}
		\displaybreak[0]\\
		&=
		\det
		\begin{bmatrix}
			A & \mathbf{a}_1 & \mathbf{a}_2\\
			\mathbf{0} & n-\ell & 0\\
			\mathbf{a}_2^* & \overline{\alpha} & n
		\end{bmatrix}
		+
		\det
		\begin{bmatrix}
			A & \mathbf{a}_1 & \mathbf{a}_2\\
			\mathbf{a}_1^* & \ell & \alpha\\
			\mathbf{a}_2^* & \overline{\alpha} & n
		\end{bmatrix}
		\displaybreak[0]\\
		&=\det \widehat{C}
		\displaybreak[0]\\
		&=\Delta_m\\
		&>0.
	\end{align*}
	By Lemma~\ref{lem:aux1}, 
	it follows that $C$ is positive definite. This implies $C\in\mathcal{K}_m$, but then $\det C>\Delta_m$ is absurd.
\end{proof}

\begin{corollary}[{\cite[Satz 2.2.]{Ehlich64b}}] Let $\widehat{C}\in \mathcal{C}_m$ be such that $\det \widehat{C}=\max\{\det C: C\in\mathcal{C}_m\}$, with $\mathcal{C}_m$ given in \eqref{eq:Cm}. Then the off-diagonal entries of $\widehat{C}$ belong to the set $\{-1,3\}$.
\end{corollary}
\begin{corollary} Let $\widehat{C}\in \mathcal{G}_m$ be such that $\det \widehat{C}=\max\{\det C: C\in\mathcal{G}_m\}$, with $\mathcal{G}_m$  given in \eqref{eq:Gm}.
	Then the off-diagonal entries of $\widehat{C}$ belong to the set $\{-1,2\}$.
\end{corollary}

\begin{definition} Let $\ell>1$ and $n\geq 2$ be integers. Define
	\[B_t = (n-\ell)I_t +\ell J_t = \begin{bmatrix}
		n & \ell & \dots & \ell\\
		\ell & n & \dots & \ell\\
		\vdots & \vdots &\ddots & \vdots\\
		\ell & \ell & \dots & n\\
	\end{bmatrix}.
	\]
	An Ehlich-type block matrix is a block matrix $M$ with blocks $B_{r_1},B_{r_2},\dots,B_{r_s}$ along its diagonal, where $r_1+r_2+\dots+r_s= n$ and whose all other entries are $-1$, i.e. we can write
	\[M=
	\begin{bmatrix}
		B_{r_1} & -1 & \dots & -1\\
		-1 & B_{r_2} & \dots & -1\\
		\vdots & \vdots &\ddots & \vdots\\
		-1 & -1 & \dots & B_{r_s}
	\end{bmatrix}.
	\]
\end{definition}
We denote the Ehlich-type block matrix with block sizes $r_1,\dots,r_s$ as $C[r_1,\dots,r_s]$.

\begin{lemma}\label{lem:aux3} Let $A$ be a Hermitian positive definite matrix of the form
	\begin{equation*}
		A=\begin{bmatrix}
			B & \mathbf{v}\\
			\mathbf{v}^* & a
		\end{bmatrix},
	\end{equation*}
	then
	\[\det A \leq a\det B,\]
	with equality if and only if $\mathbf{v}=\mathbf{0}$.
\end{lemma}
\begin{proof}
	Since $A$ is Hermitian positive definite, so is $B$. By Schur's complement formula \cite[(0.8.5.1)]{HJ90}, we have that
	\begin{equation*}
		\begin{split}
			\det A &= \det(a-\mathbf{v}^*B^{-1}\mathbf{v})\det B\\
			&=  a\det B-\mathbf{v}^*B^{-1}\mathbf{v}\det B.
		\end{split}
	\end{equation*}
	Since $B$ is positive definite, so is $B^{-1}$ therefore $\mathbf{v}^*B^{-1}\mathbf{v}\geq 0$ with equality if and only if $\mathbf{v}=0$. Thus, $\det A\leq a\det B$ with equality if and only if $\mathbf{v}=0$.
\end{proof}

\begin{theorem}\label{thm:S23} Suppose that $\widehat{C}\in\mathcal{K}_m$ satisfies $\det \widehat{C}=\Delta_m$, then up to conjugation $\widehat{C}$ is an Ehlich-type block matrix.
\end{theorem}
\begin{proof}
	We argue by contradiction, showing that if $\widehat{C}$ is not an Ehlich-type block matrix, then we can produce a matrix in $\mathcal{K}_m$ with a larger determinant. 

	Since $\det \widehat{C}=\Delta_m$, Lemma~\ref{lem:S22} shows that the off-diagonal entries of $\widehat{C}$ are in the set $\{-1,\ell\}$. If $\widehat{C}$ is not an Ehlich-type block matrix, then after simultaneous row and column permutations it can be transformed into
	\[
		C_1=\begin{bmatrix}
			A & \mathbf{a}_1 & \mathbf{a}_2 & \mathbf{a}_3\\
			\mathbf{a}_1^* & n & -1 & \ell\\
			\mathbf{a}_2^* & -1 & n & \ell\\
			\mathbf{a}_3^* & \ell & \ell & n\\
		\end{bmatrix}.
	\]
	Denote by $C_1^{(m)}$ and $C_1^{(m-1)}$ the matrices obtained from $C_1$ by deleting the $m$-th and $(m-1)$-th row respectively. Then,
	\[
		C_1^{(m)}=
		\begin{bmatrix}
			A & \mathbf{a}_1 & \mathbf{a}_2 \\
			\mathbf{a}_1^* & n & -1 \\
			\mathbf{a}_2^* & -1 & n \\
		\end{bmatrix},
		\quad\text{and }
		C_1^{(m-1)}=
		\begin{bmatrix}
			A & \mathbf{a}_1 & \mathbf{a}_3 \\
			\mathbf{a}_1^* & n & \ell \\
			\mathbf{a}_3^* & \ell & n \\
		\end{bmatrix}.
	\]
	If $\det C_{1}^{(m)}> \det C_{1}^{(m-1)}$, then we exchange the $m$-th and $(m-1)$-th rows and columns of $C_1$ to obtain the matrix
	\[C_2=
	\begin{bmatrix}
		A & \mathbf{a}_1 & \mathbf{a}_3 & \mathbf{a}_2\\
		\mathbf{a}_1^* & n & \ell & -1\\
		\mathbf{a}_3^* & \ell & n & \ell\\
		\mathbf{a}_2^* & -1 & \ell & n\\
	\end{bmatrix},
	\]
	which satisfies $\det C_2=\det C_1=\det \widehat{C}$, $C_2^{(m)}=C_1^{(m-1)}$ and $C_2^{(m-1)}=C_1^{(m)}$. Therefore, without loss of generality, we may transform the matrix $\widehat{C}$ into $C=C_1$ or $C_2$ satisying
	\begin{equation}\label{eq:t1a1}
		\det {C}^{(m)}\leq \det {C}^{(m-1)}.
	\end{equation}
	Suppose that $C=C_1$, then (\ref{eq:t1a1}) becomes:
	\begin{equation}\label{eq:t1a2}
		\det
		\begin{bmatrix}
			A & \mathbf{a}_1 & \mathbf{a}_2\\
			\mathbf{a}_1^* & n & -1 \\
			\mathbf{a}_2^* & -1 & n \\
		\end{bmatrix}
		\leq
		\det
		\begin{bmatrix}
			A & \mathbf{a}_1 & \mathbf{a}_3 \\
			\mathbf{a}_1^* & n & \ell \\
			\mathbf{a}_3^* & \ell & n \\
		\end{bmatrix}.
	\end{equation}
	Since $\det C=\Delta_m$ we have, by the last statement of Lemma~\ref{lem:S21}, that the matrix $D$ obtained from $C$ by replacing its $(m,m)$-th entry by $\ell$ is positive definite. In this case,
	\[D=
	\begin{bmatrix}
		A & \mathbf{a}_1 & \mathbf{a}_2 & \mathbf{a}_3\\
		\mathbf{a}_1^* & n & -1 & \ell\\
		\mathbf{a}_2^* & -1 & n & \ell\\
		\mathbf{a}_3^* & \ell &\ell & \ell
	\end{bmatrix}.
	\]
	By subtracting the last row from the second-to-last row and then the last column from the second-to-last column, we find
	\begin{equation*}
		\begin{split}
			\det D &=
			\det
			\begin{bmatrix}
				A & \mathbf{a}_1 & \mathbf{a}_2-\mathbf{a}_3 & \mathbf{a}_3\\
				\mathbf{a}_1^* & n & -(\ell+1) & \ell\\
				\mathbf{a}_2^*-\mathbf{a}_3^* & -(\ell+1) & n-\ell & 0\\
				\mathbf{a}_3^* & \ell & 0 & \ell
			\end{bmatrix}\\
			&=
			\det
			\begin{bmatrix}
				A & \mathbf{a}_1 & \mathbf{a}_3 &\mathbf{a}_2-\mathbf{a}_3\\
				\mathbf{a}_1^* & n & \ell & -(\ell+1)\\
				\mathbf{a}_3^* & \ell & \ell & 0\\
				\mathbf{a}_2^*-\mathbf{a}_3^* & -(\ell+1) & 0 & n-\ell \\
			\end{bmatrix}.
		\end{split}
	\end{equation*}
	Let $D'$ be the matrix in the right-hand-side above. Since $D'$ has been obtained from $D$ by a series of hermitian row/column permutations, it is also Hermitian positive definite. Applying Lemma~\ref{lem:aux3} to $D'$ with $\mathbf{v}=[\mathbf{a}_2-\mathbf{a}_3; -(\ell+1), 0]^{\intercal}\neq 0$ we find
	\begin{equation}\label{eq:t1a3}
		\det
		\begin{bmatrix}
			A & \mathbf{a}_1 & \mathbf{a}_2 & \mathbf{a}_3\\
			\mathbf{a}_1^* & n & -1 & \ell\\
			\mathbf{a}_2^* & -1 & n & \ell\\
			\mathbf{a}_3^* & \ell &\ell & \ell
		\end{bmatrix}
		=\det D'
		< (n-\ell)
		\det
		\begin{bmatrix}
			A & \mathbf{a}_1 & \mathbf{a}_3\\
			\mathbf{a}_1^* & n & \ell\\
			\mathbf{a}_3^* &\ell & \ell
		\end{bmatrix}.
	\end{equation}
	Now consider the matrix
	\[C_0=
	\begin{bmatrix}
		A & \mathbf{a}_1 & \mathbf{a}_3 & \mathbf{a}_3\\
		\mathbf{a}_1^* & n & \ell & \ell\\
		\mathbf{a}_3^* & \ell & n & \ell\\
		\mathbf{a}_3^* & \ell & \ell & n
	\end{bmatrix}.
	\]
	We find
	\begin{equation*}
		\begin{split}
			\det C_0 &=
			\det
			\begin{bmatrix}
				A & \mathbf{a}_1 & \mathbf{a}_3 & \mathbf{a}_3\\
				\mathbf{a}_1^* & n & \ell & \ell\\
				\mathbf{a}_3^* & \ell & n & \ell\\
				\mathbf{a}_3^* & \ell & \ell & n
			\end{bmatrix}\\
			&=
			\det
			\begin{bmatrix}
				A & \mathbf{a}_1 & \mathbf{a}_3 & \mathbf{a}_3\\
				\mathbf{a}_1^* & n & \ell & \ell\\
				\mathbf{a}_3^* & \ell & n & \ell\\
				\mathbf{0} & 0 & 0 & n-\ell
			\end{bmatrix}
			+
			\det
			\begin{bmatrix}
				A & \mathbf{a}_1 & \mathbf{a}_3 & \mathbf{a}_3\\
				\mathbf{a}_1^* & n & \ell & \ell\\
				\mathbf{a}_3^* & \ell & n & \ell\\
				\mathbf{a}_3^* & \ell & \ell & \ell
			\end{bmatrix}\\
			&=
			\det
			\begin{bmatrix}
				A & \mathbf{a}_1 & \mathbf{a}_3 & \mathbf{a}_3\\
				\mathbf{a}_1^* & n & \ell & \ell\\
				\mathbf{a}_3^* & \ell & n & \ell\\
				\mathbf{0} & 0 & 0 & n-\ell
			\end{bmatrix}
			+
			\det
			\begin{bmatrix}
				A & \mathbf{a}_1 & \mathbf{a}_3 & \mathbf{a}_3\\
				\mathbf{a}_1^* & n & \ell & \ell\\
				\mathbf{0} & 0 & n-\ell & 0\\
				\mathbf{a}_3^* & \ell & \ell & \ell
			\end{bmatrix}
			+
			\det
			\begin{bmatrix}
				A & \mathbf{a}_1 & \mathbf{a}_3 & \mathbf{a}_3\\
				\mathbf{a}_1^* & n & \ell & \ell\\
				\mathbf{a}_3^* & \ell & \ell & \ell\\
				\mathbf{a}_3^* & \ell & \ell & \ell
			\end{bmatrix}\\
			&=
			(n-\ell)\det
			\begin{bmatrix}
				A & \mathbf{a}_1 & \mathbf{a}_3 \\
				\mathbf{a}_1^* & n & \ell \\
				\mathbf{a}_3^* & \ell & n 
			\end{bmatrix}
			+
			(n-\ell)\det
			\begin{bmatrix}
				A & \mathbf{a}_1 & \mathbf{a}_3 \\
				\mathbf{a}_1^* & n & \ell \\
				\mathbf{a}_3^* & \ell & \ell
			\end{bmatrix}\\
			&>
			(n-\ell)\det
			\begin{bmatrix}
				A & \mathbf{a}_1 & \mathbf{a}_2\\
				\mathbf{a}_1^* & n & -1 \\
				\mathbf{a}_2^* & -1 & n \\
			\end{bmatrix}
			+
			\det
			\begin{bmatrix}
				A & \mathbf{a}_1 & \mathbf{a}_2 & \mathbf{a}_3\\
				\mathbf{a}_1^* & n & -1 & \ell\\
				\mathbf{a}_2^* & -1 & n & \ell\\
				\mathbf{a}_3^* & \ell &\ell & \ell
			\end{bmatrix}
			\quad
			\text{by (\ref{eq:t1a2}) and (\ref{eq:t1a3})}\\
			&=
			\det
			\begin{bmatrix}
				A & \mathbf{a}_1 & \mathbf{a}_2 & \mathbf{a}_3\\
				\mathbf{a}_1^* & n & -1 & \ell\\
				\mathbf{a}_2^* & -1 & n & \ell\\
				\mathbf{a}_3^* & 0 & 0 &n-\ell
			\end{bmatrix}
			+
			\det
			\begin{bmatrix}
				A & \mathbf{a}_1 & \mathbf{a}_2 & \mathbf{a}_3\\
				\mathbf{a}_1^* & n & -1 & \ell\\
				\mathbf{a}_2^* & -1 & n & \ell\\
				\mathbf{a}_3^* & \ell &\ell & \ell
			\end{bmatrix}\\
			&=
			\det C_1\\
			&=\Delta_m\\
			&>0.
		\end{split}
	\end{equation*}
	Therefore, since $C_0^{(m)}=C^{(m-1)}$ is positive definite and $\det C_0>0$, Lemma \ref{lem:aux1} implies $C_0$ is positive definite. Therefore, $C_0\in\mathcal{K}_m$, but then $\det C_0>\Delta_m$ is absurd. The case where $C=C_2$ is analogous.
\end{proof}

\begin{lemma}\label{lem:aux4}
	Let $m\geq 1$ be an integer, and let $r_1,\dots,r_s$ be positive integers with $r_1+\dots+r_s=m$. The determinant of the matrix $D=(a_{ij}J_{r_i\times r_j})_{ij}+tI_m$ is
	\begin{equation}
		t^{m-s}\det(B+tI_s),
	\end{equation}
	where $B_{ij}=r_i a_{ij}$.
\end{lemma}
\begin{proof}
	Let $W$ be the vector space spanned by vectors which are constant on the blocks of $D$. Let $m_0=0$, and $m_j=r_1+\dots+r_j$ for $j\geq 1$. 
	Then a specific basis for $W$ is given by the vectors $v_j$, $1\leq j\leq s$, which are $1$ in their entries indexed from $(m_{j-1}+1)$ to $m_j$ and $0$ everywhere else. The matrix $(a_{ij}J_{r_i\times r_j})_{ij}$ acts as $B=(r_ia_{ij})_{ij}$ on $W$ and acts as the zero matrix on $W^{\perp}$. On the other hand, $tI_m$ acts as $tI_{s}$ on $W$ and as $tI_{m-s}$ on $W^{\perp}$. This observation shows that we can conjugate $D$ by a change of basis matrix so that the resulting matrix has a $2\times 2$ block form,
	\[
		\begin{bmatrix}
			B + tI_s & 0\\
			0 & tI_{m-s}
		\end{bmatrix}.
	\]
	Therefore,
	\[\det D = \det\begin{bmatrix}
		B + tI_s & 0\\
		0 & tI_{m-s}
	\end{bmatrix}
	=t^{m-s}\det(B+tI_s),\]
	as we wanted to show.
\end{proof}
\begin{lemma} Let $C$ be an $m\times m$ Ehlich-type block matrix with $s$ blocks $B_{r_1},B_{r_2},\dots,B_{r_s}$ along the diagonal, i.e. $C=C[r_1,\dots,r_s]$ where $r_1+r_2+\dots+r_s=m$. Then,
	\begin{equation}\label{eq:det}
		\det C = (n-\ell)^{m-s}
		\left(1-\sum_{i=1}^s\frac{r_i}{(n-\ell)+(\ell+1)r_i}\right)
		\prod_{i=1}^s(n-\ell+(\ell+1)r_i).
	\end{equation}
\end{lemma}
\begin{proof} Apply Lemma~\ref{lem:aux4} with $a_{ii}=\ell$, $a_{ij}=-1$ for $i\neq j$, and $t=n-\ell$. Then,
	\[
		\det C = (n-\ell)^{m-s}\det(B+tI_s),
	\]
	where $B=(r_ia_{ij})_{ij}$. The determinant of $B+tI_s$ can be calculated as follows
	\begin{equation*}
		\begin{split}
			&\det (B+tI_s)\\
			&=
			\det
			\begin{bmatrix}
				(n-\ell)+\ell r_1 & -r_1 & \dots & -r_1\\
				-r_2 & (n-\ell)+\ell r_2 & \dots & -r_2\\
				\vdots & \vdots & \ddots & \vdots\\
				-r_s & -r_s & \dots &  (n-\ell)+\ell r_s\\
			\end{bmatrix}\\
			&=
			\det
			\begin{bmatrix}
				1 & 1 & 1 & \dots & 1\\
				0 & (n-\ell)+\ell r_1 & -r_1 & \dots & -r_1\\
				0 &-r_2 & (n-\ell)+\ell r_2 & \dots & -r_2\\
				\vdots & \vdots & \vdots & \ddots & \vdots\\
				0 &-r_s & - r_s & \dots &  (n-\ell)+\ell r_s\\
			\end{bmatrix}\\
			&=
			\det
			\begin{bmatrix}
				1 & 1 & 1 & \dots & 1\\
				r_1 & (n-\ell)+(\ell+1) r_1 & 0 & \dots & 0\\
				r_2 & 0 & (n-\ell)+(\ell+1) r_2 & \dots & 0\\
				\vdots  & \vdots & \vdots & \ddots & \vdots\\
				r_s & 0 & 0 & \dots &  (n-\ell)+(\ell+1) r_s\\
			\end{bmatrix}\\
			&=
			\det
			\begin{bmatrix}
				1-\sum_{i=1}^s\frac{r_i}{(n-\ell)+(\ell+1)r_i} & 1 
				& \dots & 1\\
				0 & (n-\ell)+(\ell+1) r_1 
				& \dots & 0\\
				\vdots  & \vdots 
				& \ddots & \vdots\\
				0 & 0 
				& \dots &  (n-\ell)+(\ell+1) r_s\\
			\end{bmatrix}\\
			&=\left(1-\sum_{i=1}^s\frac{r_i}{(n-\ell)+(\ell+1)r_i}\right)\prod_{i=1}^s((n-\ell)+(\ell+1)r_i).
		\end{split}
	\end{equation*}
	This shows Equation (\ref{eq:det}).
\end{proof}

Letting $\ell=3$ we recover Satz 3.1. of \cite{Ehlich64b}, namely that
\[\det C = (n-3)^{m-s}\left(1-\sum_{i=1}^s \frac{r_i}{n-3+4r_i}\right)\prod_{i=1}^s(n-3+4r_i),\]
for an Ehlich-type block matrix in the set $\mathcal{C}_m$ which includes Gram matrices of $\pm 1$ matrices of order $n\equiv 3\pmod{4}$.\\

The analogous determinant for an Ehlich-type block matrix in $\mathcal{G}_m$, which includes Gram matrices of matrices in $\Mat_n(\mu_3)$ with
$n\equiv 2\pmod{3}$, is obtained by letting $\ell=2$. We find:
\begin{equation}\label{eq:det3}
	\det C = (n-2)^{m-s}\left(1-\sum_{i=1}^s\frac{r_i}{n-2+3r_i}\right)\prod_{i=1}^s(n-2+3r_i).
\end{equation}

Until now, the integer parameter $n$ was used to specify the diagonal entries of a matrix in $K_m$. We now specialize $m=n$, that is, we consider $n\times n$ matrices with diagonal entries $n$. We will show that the maximal determinant in $\mathcal{K}_n$ is attained by Ehlich-type block matrices with block sizes $r$ or $r+1$. The determinant of $C[r_1,\dots,r_s]$ where $\{r_1,\dots,r_s\}\subseteq\{r,r+1\}$ will be denoted $D_n(s)$. For a fixed order $n$, this determinant depends only on the number of blocks $s$: without loss of generality, let $u\geq 1$ be the number of blocks of size $r$ and $v\geq 0$ the number of blocks of size $r+1$, then clearly $s=u+v$ and
\[ur+v(r+1) = sr+v = n.\]
From here, it follows that $n \equiv v\pmod{s}$, so, given $s$, $r$ is uniquely determined as
\[r = \left\lfloor \frac{n}{s}\right\rfloor,\]
and now $v$ and $u$ are uniquely determined as $v=n-sr$ and $u=s-v$.
%
\begin{theorem}\label{thm:S32} Let $\Delta_n=\max\{\det C:C\in\mathcal{K}_n\}$, we have
	\[\Delta_n = \max\{D_n(s):1\leq s\leq n\},\]
	where
    \begin{equation}\label{eq:Dn1}
        D_n(1)=(n-\ell)^{n-1}(n-\ell+\ell n),
    \end{equation}
    and
	\begin{equation}\label{eq:Dns}
		\begin{split}
			D_n(s)=&(n-\ell)^{n-s}(n-\ell+(\ell+1)r)^u(n+1+(\ell+1)r)^v\\
			&\times\left(1-\frac{ur}{n-\ell+(\ell+1)r}-\frac{v(r+1)}{n+1+(\ell+1)r}\right),
		\end{split}
	\end{equation}
	for $s>1$, where $r=\lfloor n/s\rfloor$, $u=s(r+1)-n$, $v=s-u$, $u\geq 1$, and $v\geq 0$. 
\end{theorem}
\begin{proof}
    By \eqref{eq:det3}, the function $D_n(s)$ is the determinant of an Ehlich-type block matrix in $\mathcal{G}_n$ with one or two block sizes, namely $\{r_1,\dots,r_s\}=\{r\}$ or $\{r,r+1\}$. Hence, by definition, $\Delta_n\geq \max\{D_n(s):1\leq s\leq n\}$. To prove the converse, let $\widehat{C}\in\mathcal{K}_n$ be such that $\det\widehat{C}=\Delta_n$. By Theorem \ref{thm:S23}, we may assume that $\widehat{C}$ is an Ehlich-type block matrix. Let $s$ be the number of blocks in $\widehat{C}$, and $\{r_1,\dots, r_s\}$ its block sizes. If $s=1$, then $\det\widehat{C}=D_n(1)$ and we are done. So, let us assume that $s>1$; we will show that all the $r_i$ are either $r$ or $r+1$ for some $r$.

	 Without loss of generality, we may assume that $r_1$ is the smallest block size in $\widehat{C}$, i.e. $r_1=\min\{r_i: 1\leq i \leq s\}$, and that $r_2$ is the largest block size in $\widehat{C}$, i.e. $r_2=\max\{r_i: 1\leq i\leq s\}$. Let $r=r_1+r_2$, we will consider the determinant in (\ref{eq:det}) as a function of $x$ where $r_1=x$, $r_2=r-x$, maintaining $r$ as a constant. We have,
	\begin{equation*}
		\begin{split}
			\det \widehat{C} = &\alpha(n-\ell)^{n-s}(n-\ell+(\ell+1)x)(n-\ell+(\ell+1)(r-x))\\
			&\times\left(1-\beta-\frac{x}{n-\ell+(\ell+1)x}-\frac{r-x}{n-\ell+(\ell+1)(r-x)}\right)
		\end{split}
	\end{equation*}
	where $\alpha=\prod_{i=3}^s(n-\ell+(\ell+1)r_i)$ and 

	\begin{equation}\label{eq:beta}
		\beta=\sum_{i=3}^s\frac{r_i}{n-\ell+(\ell+1)r_i}    
	\end{equation}

	are constants with respect to $x$. We consider the function
	\begin{equation}\label{eq:f}
		\begin{split}
			f(x)=&(n-\ell+(\ell+1)x)(n-\ell+(\ell+1)(r-x))\\
			&\times \left(1-\beta-\frac{x}{n-\ell+(\ell+1)x}-\frac{r-x}{n-\ell+(\ell+1)(r-x)}\right).
		\end{split}
	\end{equation}
	Since $r_1 \leq r_2$, 
	we have $\det \widehat{C}=\alpha f(x)$ for some $x\in[0,r/2]$.
	So, we take the domain of $f(x)$ to be the interval $[0,r/2]$.
	The function $f$ is a quadratic polynomial in $x$ and its derivative is
	\begin{equation}
		f'(x)= (\ell+1)(r-2x)((\ell-1)-(\ell+1)\beta),
	\end{equation}
	so the second derivative of $f$ is
	\[f''(x)=-2(\ell+1)((\ell-1)-(\ell+1)\beta).\]
	Since $\ell+1>0$, the point $x=r/2$ is a minimum if $(\ell-1)-(\ell+1)\beta <0$ and a maximum if $(\ell-1)-(\ell+1)\beta>0$.

	If $(\ell-1)-(\ell+1)\beta< 0$, then $x=r/2$ is a global minimum for $f$. In this case, the maximum will be achieved by $x=0$, so the determinant of $C$ may be increased by considering $s-1$ blocks instead of $s$. If $(\ell-1)-(\ell+1)\beta=0$, then $f$ is constant on $x$, and again we may choose $x=0$ so that the determinant of $C$ is not decreased by considering $s-1$ blocks instead of $s$.

	If $(\ell-1)-(\ell+1)\beta>0$, then $x=r/2$ is a global maximum, and we find that
	\[x=r_1=\left\lfloor \frac{r}{2}\right\rfloor,\quad r_2=\left\lceil \frac{r}{2}\right\rceil=\left\lfloor\frac{r+1}{2}\right\rfloor,\]
	in which case, $r_2=r_1+1$ and all other block sizes belong to $\{r_1,r_2\}$.

	We may iteratively apply the argument above until $(\ell-1)-(\ell+1)\beta>0$ in which case we arrive at a situation where there are only two block sizes differing by $1$, or until $s=1$ in which case there is a unique block of size $n$. In either case, we find that the maximal determinant is achieved by an Ehlich-type block matrix with $u\geq 1$ blocks of size $r$ and $v\geq 0$ blocks of size $r+1$, where $u+v=s$ and $ur+v(r+1)=n$. The associated determinant is $D_n(s)$ given in (\ref{eq:Dns}). 
\end{proof}

With this, the proof of Theorem \ref{thm:main} consists in narrowing down the possible values of $s$ which maximize the determinant of $D_n(s)$ for all $n$. This analysis is dependent on the specific value of $\ell$, so from now on we consider $\ell=2$, the case $\ell=3$ is analogous and has been developed in \cite{Ehlich64b}. When $\ell=2$, we have that Equations \eqref{eq:Dn1} and \eqref{eq:Dns} become
\begin{equation}\label{eq:Dns2}\hspace{-2pt}
    D_n(s)=\begin{cases}
        (n-2)^{n-1}(3n-2) & \text{if } s=1,\\
        (n-2)^{n-s}(n-2+3r)^u(n+1+3r)^v\left(1-\frac{ur}{n-2+3r}-\frac{v(r+1)}{n+1+3r}\right) &\text{if } s>1.
    \end{cases}
\end{equation}
\begin{lemma}\label{lem:3A}
	Let $\Delta_n=\max\{\det C: C\in \mathcal{G}_n\}$. For all integers $n\geq 5$ with $n\equiv 2\pmod{3}$, we have
    \[\Delta_n=\max\{D_n(s): 1 \leq s \leq n\}=\max\{D_n(3),D_n(4)\}.\]
\end{lemma}
\begin{proof}
    Let $\widehat{C}\in\mathcal{G}_n$ be a matrix satisfying $\det \widehat{C}=\Delta_n$, then by Theorem \ref{thm:S32} with $\ell=2$, we have 
    \[\det \widehat{C}=D_n(s),\]
     where $D_n(s)$ is given by \eqref{eq:Dns2}, for some integer $s$ with $1\leq s \leq n$. First, suppose that $s\geq 5$ and introduce a new variable in the function $D_n(s)$ as follows:
    \[D_n(s,x)=(n-2)^{n-s-1}(n-2+3r)^{u-1}(n+1+3r)^vg(x),\]
    where
    \[g(x)=(n-2+3x)(n-2+3r-3x)\left(1-\beta-\frac{x}{n-2+3x}-\frac{r-x}{n-2+3r-3x}\right),\]
and
    \[\beta=\begin{cases}
        \frac{u-1}{n-2+3r} & \text{if } v=0,\\
        \frac{(u-1)r}{n-2+3r}+\frac{(v-1)(r+1)}{n+1+3r} & \text{if } v>0.
    \end{cases}
    \]
    For an integer value of $x$, the value $D_n(s,x)$ is equal to the determinant of the Ehlich-type block matrix $C'\in\mathcal{G}_n$ obtained from $\widehat{C}$ by replacing one of its blocks of size $r$ with two blocks of size $x$ and $r-x$. Then, by maximality of $\widehat{C}$, 
    \[\det C'=D_n(s,x)\leq D_n(s,0)=D_n(s)=\det\widehat{C}=\Delta_n.\]
     Hence, since $g$ is an at most quadratic polynomial with $g(0)=g(r)$ we must have $g''(x)\leq 0$ for all $x\in (0,r)$. The second derivative of $g(x)$ is
    \[g''(x)=-6(1-3\beta).\]
    Hence, the maximality of $\widehat{C}$ implies $1-3\beta \geq 0$.

    If $s\geq 5$, then
    	\[r\leq \frac{n}{s}\leq \frac{n}{5}<\frac{n-1}{3}\quad \text{since } n\geq 5.\]
    Then, if $v=0$
	\[(n-2+3r)(1-3\beta)=-2n-2+6r< 0.\]
    If $v>0$, then
    \begin{align*}
	    (n+1+3r)(1-3\beta)&\leq n+1+3r -3(u-1)r-3(v-1)(r+1)\\
	    &=-2n+4+9r< 0.
    \end{align*}
    So, in either case, $s\geq 5$ implies $1-3\beta <0$. This implies $s\leq 4$.
    
    To show $s\ne 2$, observe that if $s= 2$, then $sr=2r=n-v$ and $v\leq s-1$ imply
    \[r= \frac{n-v}{2}\geq \frac{n}{2}>\frac{n+1}{3},\quad\text{for } n>2.\]
    Notice that for $v\geq 0$
    \[\beta\leq \frac{(u-1)r}{n-2+3r}+\frac{v(r+1)}{n+1+3r}.\]
    Therefore,
    \begin{align*}
    (n-2+3r)(1-3\beta) &\geq (n-2+3r)\left(1-3\left(\frac{(u-1)r}{n-2+3r}+\frac{v(r+1)}{n-2+3r}\right)\right)\\
    &=n-2+3r-3(ur+v(r+1)-r)\\
    &=n-2+3r-3(n-r)\\
    &=-2n-2+6r>0.
    \end{align*}
    Which contradicts $1-3\beta\geq 0$.

    Finally, to show that $s\neq 1$, consider the function
    \[D(x)=(n-2)^{n-2}(n-2+3x)(n-2+3(n-x))(1-\frac{x}{n-2+3x}-\frac{n-x}{n-2+3(n-x)}),\]
    where $0\leq x \leq \lfloor n/2\rfloor$. Notice that $D(0)=D_n(1)$ and $D(\lfloor n/2 \rfloor)=D_n(2)$. 
    Write $D(x)$ as
    \[D(x)=(n-2)^{n-2}h(x),\]
    where
    \[h(x)=3x(n-x)+(3n^2-8n+4).\]
    Then,
    \[h'(x)=3(n-2x),\]
    and $h''(x)=-6<0$, which implies that the maximum of $D(x)$ is achieved at $x=\lfloor n/2\rfloor$, which in turn implies that $D_n(2)>D_n(1)$.

    This shows that $s\geq 3$.
\end{proof}
\begin{lemma}\label{lem:3B}
    For $n\equiv 2\pmod{3}$ with $n\geq 26$, 
    \[D_n(3)<D_n(4).\]
	For $n\equiv 2\pmod{3}$ with $5\leq n\leq 23$,
	\[D_n(3)>D_n(4).\]
\end{lemma}
\begin{proof}
    From $ur+v(r+1)=sr+v=n$, we have that
    \[1-\frac{n}{n-2+3r}\leq 1-\frac{ur}{n-2+3r}-\frac{v(r+1)}{n+1+3r}< 1-\frac{n}{n+1+3r}.\]
    Combining the above with \eqref{eq:Dns2},
    \begin{align}(n-2)^{-n+s}D_n(s)&\geq (n-2+3r)^u(n+1+3r)^v\left(1-\frac{n}{n-2+3r}\right)\nonumber\\
    &=(n-2+3r)^{u-1}(n+1+3r)^v(3r-2).\label{eq:lb}
    \end{align}
    Additionally,
    \begin{equation}\label{eq:ub}
        (n-2)^{-n+s}D_n(s)\leq (n-2+3r)^u(n+1+3r)^{v-1}(3r+1).
    \end{equation}

    Let $s=4$ in \eqref{eq:lb} and multiply by $4^4$, then using the fact that $4r=sr=n-v$ and $u=4-v$, we have
    \begin{align*}
        4^4(n-2)^{-n+4}D_n(4) &\geq (4n-8+3(n-v))^{3-v}(4n+4+3(n-v))^v(3(n-v)-8)\\
        &=(7n-8-3v)^{3-v}(7n+4-3v)^v(3n-8-3v).
    \end{align*}
    Denote the right-hand-side as $A_v=(7n-8-3v)^{3-v}(7n+4-3v)^v(3n-8-3v)$, then
\begin{align*}
&A_v-A_3
\\&=(7n-8-3v)^{3-v}(7n+4-3v)^v(3n-8-3v)-(7n-5)^3(3n-17)
\\&=
\begin{cases}
1764n^3 - 630n^2 - 2988n + 1971
&\text{if $v=0$,}\\
1176n^3 - 672n^2 + 2040n - 3456
&\text{if $v=1$,}\\
588n^3 - 462n^2 + 3252n - 1341
&\text{if $v=2$}
\end{cases}
\\&\geq0
\end{align*}
for $n\geq2$.
Thus,
\begin{align*}
4^4(n-2)^{-n+4}D_n(4) &\geq 
(7n-5)^3(3n-17)\\
&= 7^3\cdot 3\left(n-\frac{5}{7}\right)^3\left(n-\frac{17}{3}\right)\\
&> 7^3\cdot 3(n-1)^3(n-6).
\end{align*}       
    
    Now, let $s=3$ in \eqref{eq:ub}, since $3r=n-v$ and $u=3-v$,
    \[(n-2)^{-n+3}D_n(3)< (2n-2-v)^{3-v}(2n+1-v)^{v-1}(n+1-v).\]
    Let $B_v=(2n-2-v)^{3-v}(2n+1-v)^{v-1}(n+1-v)$, then
    \begin{align*}
    (2n+1)(B_0-B_v)&=2^3(n-1)^3(n+1)-(2n-2-v)^{3-v}(2n+1-v)^{v-1}(n+1-v)(2n+1)\\
    &=\begin{cases}
    4n^3-6n^2+7n-8 & \text{ if } v=1\\
    8n^3-14n^2+10n-4 &\text{ if } v=2
    \end{cases}\\
    &\geq 0,
    \end{align*}
    for $n\geq 2$, so $B_v\leq B_0$ for $n\geq 2$. This implies,
    \begin{align*}
        (n-2)^{-n+3}D_n(3)&\leq 2^3(n-1)^3(n+1)(2n+1)^{-1}\\
        &=2^2(n-1)^3(n+1)\left(n+\frac{1}{2}\right)^{-1}.
    \end{align*}
    
    To conclude, we use the bounds on $D_n(4)$ and $D_n(3)$ to obtain
    \begin{align*}
        \frac{D_n(4)}{D_n(3)}&>\frac{7^3\cdot 3}{4^5}\frac{(n-1)^3(n-6)(n+1/2)}{(n-2)(n-1)^3(n+1)}\\
        &=\frac{7^3\cdot 3}{4^5}\frac{(n-6)(n+1/2)}{(n-2)(n+1)}\\
        &=\frac{7^3\cdot 3}{4^5}\left(1-\frac{4}{n-2}\right)\left(1-\frac{1/2}{n+1}\right)\\
        &>1,
    \end{align*}
    for $n\geq 928$. There are 308 integers $n$ with $5\leq n < 928$ and $n\equiv 2\pmod{3}$, for each such value of $n$ we compared $D_n(3)$ and $D_n(4)$ by direct computation and found that $D_n(4)>D_n(3)$ for $n\geq 26$ and $D_n(4) < D_n(3)$ for $5\leq n \leq 23$. 
\end{proof}

\begin{proof}[Proof of Theorem \ref{thm:main}] Let $M\in \Mat_n(\mu_3)$, then by Proposition~\ref{prop:1}, there exists a diagonal matrix $P\in\Mat_n(\mu_3)$ such that $(PM)(PM)^*\in\mathcal{G}_n$. In particular, 
	\begin{equation}\label{eq:p1}
		|\det M|^2\leq \Delta_n.
	\end{equation}
	By Theorem~\ref{thm:S32}, we have that $\Delta_n=\max\{D_n(s):1\leq s\leq n\}$. Then, if $n=2$,
	\begin{equation}\label{eq:p2}
		\Delta_n=\max\{D_n(1),D_n(2)\}=D_n(2),
	\end{equation}
	since $D_2(1)=0<3=D_2(2)$. If $n\geq 5$, then by Lemma~\ref{lem:3A}, $\Delta_n=\max\{D_n(3),D_n(4)\}$, and by Lemma~\ref{lem:3B} it follows that
	\begin{equation}\label{eq:p3}
		\Delta_n = 
		\begin{cases}
			D_n(3) & \text{ if } 5\leq n \leq 23,\\
			D_n(4) & \text{ if } n\geq 26.
		\end{cases}
	\end{equation}
	Combining \eqref{eq:p1} with \eqref{eq:p2} and \eqref{eq:p3} we obtain the desired result. \qedhere
\end{proof}

\section{Asymptotics}
In this section, we use the Landau big $O$ notation, where a real-valued function $f$ is said to be $O(g)$, written $f=O(g)$, if there exists a natural number $N$ and a constant $b\in \R$ such that for all $x\geq N$, $|f(x)|\leq bg(x)$; here $g$ is a real-valued function satisfying $g(x)>0$ for all $x$.

Recall, from Equation~\eqref{eq:Dns2}, that
\[D_n(4)=(n-2)^{n-4}(n-2+3r)^u(n+1+3r)^v\left(1-\frac{ur}{n-2+3r}-\frac{v(r+1)}{n+1+3r}\right),\]
where $r=\lfloor n/4\rfloor$, and $u+v=4$. Then, $r=n/4+O(1)$, so $n-2+3r=7n/4+O(1)$ and $n+1+3r=7n/4+O(1)$. Therefore,
\begin{align*}
    D_n(4)&=(n-2)^{n-4}(7n/4+O(1))^{u+v}\left(1-\frac{un/4+O(1)}{7n/4+O(1)}-\frac{vn/4+O(1)}{7n/4+O(1)}\right)\\
    &=(n-2)^{n-4}(7n/4+O(1))^4\left(1-(u+v)/7+O(1/n)\right)\\
    &=(7/4)^4(n-2)^{n-4}(n+O(1))^4(3/7+O(1/n))\\
    &=(7/4)^4(3/7)(n-2)^{n-4}(n+O(1))^4(1+O(1/q)\\
    &=\frac{1029}{256}(n-2)^{n-4}(n^4+O(n^3)).
\end{align*}
From here, it follows that
\begin{align*}
    \lim_{n\rightarrow\infty} \frac{D_n(4)}{n^n} &= \frac{1029}{256}\lim_{n\rightarrow\infty}\frac{(n-2)^{n-4}(n^4+O(n^3))}{n^n}\\
    &=\frac{1029}{256}\lim_{n\rightarrow\infty}\left(\frac{n-2}{n}\right)^n\frac{n^4+O(n^3)}{(n-2)^4}\\
    &=\frac{1029}{256}\lim_{n\rightarrow\infty}\left(1-\frac{2}{n}\right)^n\frac{n^4+O(n^3)}{n^4+O(n^3)}\\
    &=\frac{1029}{256e^2}\approx 0.5440.
\end{align*}
Taking square-roots we find that the asymptotic ratio of the determinant bound of Theorem~\ref{thm:main} to the Hadamard bound is $\sqrt{1029}/16e\approx 0.7375$. In summary, the sharpest known upper bounds for the maximal determinant $\gamma_3(n)$ of a matrix of order $n$ with entries in $\mu_3$ are 
\[
\gamma_3(n)\leq
\begin{cases}
    n^{n/2} & \text{if } n\equiv 0\pmod 3,\\
    \sqrt{2n-1}(n-1)^{(n-1)/2} & \text{if } n\equiv 1 \pmod 3,\\
    \sqrt{D_n(4)} &\text{if } n\equiv 2\pmod 3.
\end{cases}
\]
The asymptotic ratios to the Hadamard bound when $n\rightarrow\infty$ are $1$, $\sqrt{2/e}\approx0.8578$ and $\sqrt{1029}/16e\approx 0.7375$, for $n\equiv 0,1,$ and $2\pmod 3$, respectively. 

In \cite[Corollary 6]{Ponasso25}, a lower bound for the maximal determinant of a matrix of order $q+1\equiv 2\pmod{3}$, where $q$ is a prime power, was given. For a prime power $q\equiv 1\pmod{3}$, let $Q_q$ be the cubic generalized Paley core; this is the matrix indexed by elements of $\F_q$, with $(i,j)$-th entry equal to:
\[(Q_{q})_{ij}=\chi(i-j),\]
where $\chi$ is a primitive cubic character in $\F_q$. Define the matrix $M_q$ as follows
\[
M_q = \begin{bmatrix}
    \omega & \nf 1^{\intercal}\\
    \nf{1} & Q_q+\omega I_{q}.
\end{bmatrix}
\]
\begin{proposition}[cf. \cite{Ponasso25}, Corollary 6]\label{prop:Mq}
    Let $M_q$ be defined as above, then
    \[|\det M_q|^2 = (q^2+q+1)(q^3+qc+1)^{(q-1)/3},\]
    where $c=O(\sqrt{q})$ is an integer uniquely determined by the cubic cyclotomic numbers of $\F_q$.
\end{proposition}
Note that $(q^3+cq+1)^{(q-1)/3}=q^{q-1}(1+(cq+1)/q^3)^{(q-1)/3}$. On the other hand, since $c=O(\sqrt{q})$
\begin{align*}
\lim_{q\rightarrow\infty}q\log(1+(cq+1)/q^3) &= \lim_{q\rightarrow\infty}\frac{\log(1+(cq+1)/q^3)}{1/q}\\
    &=\lim_{q\rightarrow\infty}-\frac{cq^3-3q^2(cq+1)}{q^4+q(cq+1)}\\
    &=0.
\end{align*}
Therefore,
\[\lim_{q\rightarrow\infty}\frac{(q^3+cq+1)^{(q-1)/3}}{q^{q-1}}=\lim_{q\rightarrow\infty}(1+(cq+1)/q^3)^{q/3}(1+(cq+1)/q^3)^{-1/3}=1.\]
Hence, we can calculate the asymptotic ratio of $|\det M_q|^2$ to the Hadamard bound as
\begin{align*}
    \lim_{q\rightarrow\infty} \frac{|\det M_q|^2}{(q+1)^{q+1}} &= \lim_{q\rightarrow\infty} \frac{(q^2+q+1)(q^3+qc+1)^{(q-1)/3}}{(q+1)^{q+1}}\\
    &=\lim_{q\rightarrow\infty}\frac{(q^2+q+1)q^{q-1}}{(q+1)^{q+1}}\\
    &=\lim_{q\rightarrow\infty}\left(\frac{q}{q+1}\right)^{q}\frac{q^2+q+1}{q(q+1)}\\
    &=\lim_{q\rightarrow\infty}\left(1+\frac{1}{q}\right)^{-q}\\
    &=1/e.
\end{align*}
From this, we have
\begin{align*}
    \lim_{q\rightarrow\infty}\frac{|\det M_q|^2}{D_{q+1}(4)} &= \lim_{q\rightarrow\infty} \frac{|\det M_q|^2}{(q+1)^{q+1}}\frac{(q+1)^{q+1}}{D_{q+1}(4)}\\
    &=\frac{1}{e}\frac{256e^2}{1029}\\
    &=\frac{256e}{1029}\approx 0.6763.
\end{align*}
Therefore, the construction of Proposition~\ref{prop:Mq} achieves an asymptotic ratio of $16\sqrt{e/1029}\approx 0.8223$ to the determinant bound of Theorem~\ref{thm:main}.

The status of the maximal determinant problem for matrices over the third roots, with the bounds updated in view of Theorem \ref{thm:main} can be found in Table \ref{tab:1}. The smallest open cases in the determination of the maximal determinant over $\mu_3$ are $n=14$ and $n=15$. The largest known determinant at order $n=14$ is given by the matrix
\[
M_{14}=M_{14}^{\intercal}=\left[
\begin{array}{*{6}{c}|*{5}c|*{3}{c}} 
0&0&0&1&0&0&2&2&0&1&1&1&1&0\\
0&0&1&0&0&0&2&1&1&0&2&1&0&1\\
0&1&0&0&0&0&1&2&1&2&0&0&1&1\\
1&0&0&2&0&0&1&1&2&2&2&1&1&2\\
0&0&0&0&1&2&1&2&2&1&2&1&2&1\\
0&0&0&0&2&1&2&1&2&2&1&2&1&1\\
\hline
2&2&1&1&1&2&0&1&0&2&0&1&1&1\\
2&1&2&1&2&1&1&0&0&0&2&1&1&1\\
0&1&1&2&2&2&0&0&2&0&0&1&1&0\\
1&0&2&2&1&2&2&0&0&0&0&0&1&1\\
1&2&0&2&2&1&0&2&0&0&0&1&0&1\\
\hline
1&1&0&1&1&2&1&1&1&0&1&2&0&0\\
1&0&1&1&2&1&1&1&1&1&0&0&2&0\\
0&1&1&2&1&1&1&1&0&1&1&0&0&2
\end{array}\right].
\]
This matrix is balanced and normal, i.e. $M_{14}M_{14}^*=M_{14}^*M_{14}$. Its Gram matrix is of Ehlich block type with shape given by the parititon $(6,5,3)$ and its determinant is $|\det M_{14}|^2= 2^{24}\cdot 3^{13}\cdot 223< 2^{24}\cdot 3^{15}\cdot 5^2= D_{14}(3)$. The largest known determinant at order $15$ is given by the bordered circulant matrix
\[
M_{15}=
\left[
\begin{array}{c|*{14}{c}} 
0&1&1&1&1&1&1&1&1&1&1&1&1&1&1\\
\hline
1&0&1&1&0&1&0&2&2&1&1&2&1&2&0\\
1&0&0&1&1&0&1&0&2&2&1&1&2&1&2\\
1&2&0&0&1&1&0&1&0&2&2&1&1&2&1\\
1&1&2&0&0&1&1&0&1&0&2&2&1&1&2\\
1&2&1&2&0&0&1&1&0&1&0&2&2&1&1\\
1&1&2&1&2&0&0&1&1&0&1&0&2&2&1\\
1&1&1&2&1&2&0&0&1&1&0&1&0&2&2\\
1&2&1&1&2&1&2&0&0&1&1&0&1&0&2\\
1&2&2&1&1&2&1&2&0&0&1&1&0&1&0\\
1&0&2&2&1&1&2&1&2&0&0&1&1&0&1\\
1&1&0&2&2&1&1&2&1&2&0&0&1&1&0\\
1&0&1&0&2&2&1&1&2&1&2&0&0&1&1\\
1&1&0&1&0&2&2&1&1&2&1&2&0&0&1\\
1&1&1&0&1&0&2&2&1&1&2&1&2&0&0
\end{array}
\right].
\]
The first row of the circulant core can be constructed using the mapping
\begin{align*}
    f:\Z/14\Z &\rightarrow \Z/3\Z\\
    d &\mapsto ((1+\varepsilon)\chi(u)-\varepsilon) \pmod{3},
\end{align*}
where $(u,\varepsilon):=(d\mod 7,d\mod 2)$ and $\chi(u)$ is the quadratic character on $\F_7$, i.e. $\chi(u)=+1$ if $u$ is a non-zero square, $\chi(u)=-1$ if $u$ is a non-square, and $\chi(0)=0$. We have $|\det M_{15}|^2 = 2^{22}\cdot 3^{20}\cdot 19$, so $M_{15}$ achieves a ratio of approximately $0.7966$ of the Hadamard bound. 
A matrix monomially equivalent to $M_{15}$ with a $2\times 2$ block matrix core is
\[\left[
\begin{array}{ccc}
 0   & \nf 1^{\intercal} & \nf 1^{\intercal}\\
 \nf 1 &A & \omega^2A^{\intercal}\\
 \nf 1 &\omega^2A^{\intercal} & A
\end{array}
\right],
\]
where $A=I_7+\omega^2D+\omega(J_7-I_7-D)$ and
\[
D = \begin{bmatrix}
0&1&1&0&1&0&0\\
0&0&1&1&0&1&0\\
0&0&0&1&1&0&1\\
1&0&0&0&1&1&0\\
0&1&0&0&0&1&1\\
1&0&1&0&0&0&1\\
1&1&0&1&0&0&0
\end{bmatrix},
\]
is the incidence matrix of the Fano plane.
The maximality of $M_{14}$ and $M_{15}$ remains undecided.
\begin{table}
    \centering
    \begin{tabular}{llc|llc|llc}
    \hline
        $n$ & $|\det|^2/3^{n-1}$  & R & $n$ & $|\det|^2/3^{n-1}$  & R & $n$ & $|\det|^2/3^{n-1}$  & R  \\
        \hline
        &&&
        $1$ & $1$ & $1$ & 
        $2$ & $1$ & $1$\\
        $3$ & $3$ & $1$ &
        $4$ & $7$ & $1$ &
        $5$ & $3\times 7$ & $1$\\
        $6$ & $2^6\times 3$ & $1$ &
        $7$ & $2^6\times 13$ & $1$ &
        $8$ & $2^{12}$ & $0.99$\\
        $9$ & $3^{10}$ & $1$ &
        $10$ & $3^9\times 19$ & $1$ &
        $11$ & $3^9\times 7\times 19$ & $1$\\
        $12$ & $2^{24}\times 3$ & $1$ &
        $13$ & $2^{24}\times 5^2$ & $1$ &
        $14$ & $2^{24}\times 223$\red{??} & \red{$0.99$}\\
        $15$ & $2^{22}\times 3^6\times 19$\red{??} & \red{0.79}&
        $16$ & $2^{24}\times 3^8\times 7$\red{??} & \red{0.90}&
        $17$ & $13^5\times 67^4$\red{??} & \red{0.85}\\
        $18$ & $2^{18}\times 3^{19}$ & 1 &
        $19$ & $13\times 37^2\times 342037^2$\red{??} & \red{0.74} &
        $20$ & $7^6\times 37^6\times 127$\red{??} & \red{0.88}\\
        \hline
    \end{tabular}
    \caption{Largest known determinant values for matrices of order $n$ over the third roots: columns labeled $n$ indicate the matrix order. The columns labeled $|\det|^2/3^{n-1}$ contain the absolute value squared of the largest known determinant value at order $n$, divided by $3^{n-1}$. Columns labeled R contain the ratio of the largest known determinant to the Hadamard bound at orders $n\equiv 0\pmod{3}$, the Barba bound at orders $n\equiv 1\pmod{3}$ and the bound of Theorem \ref{thm:main} at orders $n\equiv 2\pmod{3}$. The symbol \red{??} is used whenever there is no proof of maximality known for the reported determinant value.}
    \label{tab:1}
\end{table}

\subsection*{Acknowledgements}
AM was supported by the Japan Society for the Promotion of Science (JSPS) KAKENHI Grant Number 25K07095. GNP was supported by the JSPS Postdoctoral fellowship, and by the JSPS KAKENHI Grant Number 24KF0176.

\printbibliography
\end{document}